\documentclass[11pt,titlepage]{article}
\usepackage{latexsym}
\usepackage{amsfonts}
\usepackage{amsmath}
\usepackage{amssymb}
\usepackage{graphicx}
\usepackage{tikz}
\usetikzlibrary{arrows, decorations.markings}
\newcommand\blackslug{\hbox{\hskip 1pt \vrule width 4pt height 8pt depth 1.5pt
        \hskip 1pt}}
\newcommand\bbox{\hfill \quad \blackslug \medbreak}

\newtheorem{theorem}{}[section]

\newcommand{\Proof}{\noindent{\bf Proof.}\ \ }

\tikzstyle{every node}=[circle, draw, fill=black!50, inner sep=0pt, minimum width=20pt]

\tikzstyle{input}=[circle,
                                    thick,
                                    minimum size=2.2cm,
                                    draw=black!80,
                                    fill=white!20]

\tikzstyle{matrx}=[rectangle,
                                    thick,
                                    minimum size=1.8cm,
                                    draw=black!80,
                                    fill=white!20]

\tikzstyle{vecArrow} = [thick, decoration={markings,mark=at position
   1 with {\arrow[semithick]{open triangle 60}}},
   double distance=1.8pt, shorten >= 5.5pt,
   preaction = {decorate},
   postaction = {draw,line width=1.4pt, white,shorten >= 4.5pt}]
\tikzstyle{innerWhite} = [semithick, white,line width=1.4pt, shorten >= 4.5pt]

\title{Excluding pairs of tournaments}
\author{Krzysztof Choromanski\\
Columbia University\\
New York, NY
}
\date{September 7, 2014; revised \today}

\newcommand{\Keywords}{the Erd\H{o}s-Hajnal conjecture, the regularity lemma, pairs of tournaments}

\begin{document}
\maketitle
\begin{abstract}
The Erd\H{o}s-Hajnal conjecture states that for every given undirected graph $H$ there exists a constant $c(H)>0$ such that
every graph $G$ that does not contain $H$ as an induced subgraph contains a clique or a stable set of size at least
$|V(G)|^{c(H)}$. The conjecture is still open. 
Its equivalent directed version states that for every given tournament $H$ there exists a constant $c(H)>0$ such that every $H$-free tournament $T$ contains a transitive subtournament of order at least $|V(T)|^{c(H)}$.
We prove in this paper that $\{H_{1},H_{2}\}$-free tournaments $T$ contain transitive subtournaments of size at least 
$|V(T)|^{c(H_{1},H_{2})}$ for some $c(H_{1},H_{2})>0$ and several pairs of tournaments: $H_{1}$, $H_{2}$.
In particular we prove that $\{H,H^{c}\}$-freeness
implies existence of the polynomial-size transitive subtournaments for several tournaments $H$ for which the conjecture is
still open ($H^{c}$ stands for the \textit{complement of $H$}).
To the best of our knowledge these are first nontrivial results of this type.

\end{abstract}

\maketitle {\bf Keywords:} \Keywords

\section{Introduction}

The Erd\H{o}s-Hajnal conjecture states that for every given undirected graph $H$ there exists a constant $c(H)>0$ such that
every graph $G$ that does not contain $H$ as an induced subgraph contains a clique or a stable set of size at least
$|V(G)|^{c(H)}$. The conjecture is still open. 
Its equivalent directed version states that for every given tournament $H$ there exists a constant $c(H)>0$ such that every $H$-free tournament $T$ contains 
a transitive subtournament of order at least $|V(T)|^{c(H)}$.
In the undirected setting so far the conjecture was proved for some graphs on at most five vertices and the graphs
obtained from them by the so-called \textit{substitution procedure} (see \cite{pach}).
Much more is known for the directed setting. The conjecture was proved for all tournaments on at most five vertices and
so-called \textit{galaxies} (see \cite{choromanski1}). The conjecture was then proved for the family of so-called \textit{constellations}
which contains the family of galaxies (see \cite{choromanski2}) . Even more recently the conjecture was proved for more tournaments
with the use of the so-called \textit{strong EH-property} and the notion of product tournaments (see \cite{choromanski3}).
Forbidden tournaments which absence implies existence of the linear size or near-linear size transitive subsets were fully characterized in \cite{seymour} and \cite{pseudo}. 
Instead of forbidding just one undirected graph/tournament, one can state the analogous conjecture for the case where
all the graphs from the given (possibly infinite) class $\mathcal{C}$ are forbidden. In particular one can analyze the setting, where
we forbid undirected graph $H$ and its complement $H^{c}$. It is not known whether the conjecture holds even for this scenario.
Below we list some known results regarding excluding families of undirected graphs.
Chudnovsky and Zwols proved (see: \cite{chudzwols}) that excluding four-edge path and five-edge path in the complement implies existence
of polynomial-size cliques or stable sets. This results was further refined in \cite{seymourchud}.
Surprisingly, much more general result holds. Excluding an arbitrary path and an arbitrary antipath gives the conjecture. This was
very recently proved in \cite{thom}.
Several interesting results regarding excluding pairs of graphs are included in \cite{chudscottseym}.
For comparison, in the directed case there are almost no results regarding excluding pairs of tournaments and regarding the 
Erd\H{o}s-Hajnal conjecture in this setting.

We prove in this paper that $\{H_{1},H_{2}\}$-free tournaments $T$ contain transitive subtournaments of size at least 
$|V(T)|^{c(H_{1},H_{2})}$ for some $c(H_{1},H_{2})>0$ and several pairs of tournaments: $H_{1}$, $H_{2}$.
Thus we prove the Erd\H{o}s-Hajnal conjecture for these pairs: $\{H_{1},H_{2}\}$.
In particular we prove that $\{H,H^{c}\}$-freeness
implies existence of the polynomial-size transitive subtournaments for several tournaments $H$ for which the conjecture is
still open ($H^{c}$ stands for the \textit{complement of $H$}). 
To the best of our knowledge these are first nontrivial results of this type.
Before stating our results formally, we need to introduce some notation and few definitions.

All graphs in this paper are finite and simple. Let $G$ be an undirected  
graph. The vertex 
set of $G$ is denoted by $V(G)$, and the edge set by $E(G)$. We write 
$|G|$ to mean $|V(G)|$. For undirected $G$, given $X \subseteq V(G)$, we denote by  $G|X$ the
{\em  subgraph of} $G$ {\em induced by} $X$, that is the graph with vertex set $X$, 
in which $x,y \in X$ are adjacent if and only if they are adjacent in $G$.
For an undirected  graph $H$, we say that $G$ is $H$-{\em free} if no induced 
subgraph of $G$ is isomorphic to $H$.
A {\em clique} in $G$ is a subset of $V(G)$ all of whose elements are pairwise
adjacent, and a {\em stable set} in $G$ is a subset of $V(G)$ all of whose 
elements are pairwise non-adjacent. For a graph $H$ and a vertex $v \in V(H)$
we denote by $H \backslash v$ a graph obtained from $H$ by deleting $v$ and all
edges of $H$ that are: adjacent to $v$ in the undirected setting and: adjacent to or from
$v$ in the directed setting.

The Erd\H{o}s-Hajnal Conjecture is the following: 

\begin{theorem}
\label{erdos-hajnal-conjecture-undirected}
For every undirected graph $H$ there exists a constant $c(H)>0$ such that the following holds: every $H$-free graph $G$ contains a clique or a stable set
of size at least $|G|^{c(H)}$.
\end{theorem}

A version of \ref{erdos-hajnal-conjecture-undirected} 
in the directed setting was formulated in \cite{pach}. 
To state it, we need some definitions. A {\em tournament} is a directed graph $T$, where for every two vertices $u,v$ exactly one of $(u,v)$, $(v,u)$ is an edge of $T$ (that is, a directed edge). If $(u,v) \in E(T)$, we say
that $u$ is {\em adjacent to} $v$, and that $v$ is {\em adjacent from} $u$.
A tournament 
is {\em transitive} if it contains no directed cycle (equivalently, no directed
cycle of length three). Let $T$ be a tournament.  We denote its vertex set
by $V(T)$ and its edge set by $E(T)$,  and write $|T|$ for $|V(T)|$. We
refer to $|T|$ as the {\em order} of $T$.
Given $X \subseteq V(T)$, the {\em subtournament of} $T$ {\em induced by} $X$, 
denoted by $T|X$, is the tournament with vertex set $X$, such that for 
$x,y \in X$, $(x,y)$ is a directed edge of $T|X$ if and only if $(x,y) \in E(T)$.
Given a tournament $S$, we say that $T$ {\em contains} $S$ if $S$ is isomorphic
to $T|X$ for some $X \subseteq V(T)$. If $T$ does not contains $S$, we say that
$T$ is $S$-{\em free}. 
For two disjoint subsets $A,B \subseteq V(T)$ we say that $A$ is complete to $B$
if every vertex of $A$ is adjacent to every vertex of $B$.
The conjecture from \cite{pach} is the following.
 
\begin{theorem}
\label{erdos-hajnal-conjecture-directed}
For every tournament $S$ there exists a constant $c(S)>0$ such that the following holds: every $S$-free tournament $T$ contains a transitive subtournament of order at least $|T|^{c(S)}$.
\end{theorem}

The complement of the graph $H$ will be denoted as $H^{c}$. We will use this notation when $H$ is an undirected graph or a tournament.
In the tournament setting the complement of the tournament $H$ is defined as a tournament obtained from $H$ by reversing directions of all the edges.
One can easily see that if the conjecture is true for $H$ then it is true also for $H^{c}$.

For a class $\mathcal{C}$ of undirected graphs and an undirected graph $G$ we say that $G$ is \textit{$\mathcal{C}$-free} if it is $C$-free for every $C \in \mathcal{C}$.
For a class $\mathcal{C}$ of tournaments and a tournament $T$ we say that $T$ is \textit{$\mathcal{C}$-free} if it is $C$-free for every $C \in \mathcal{C}$.

When $\mathcal{C}$ consists of undirected graphs we say that $\mathcal{C}$ has the \textit{Erd\H{o}s-Hajnal property} if the following holds: 
there exists $c_{\mathcal{C}}>0$ such that every $\mathcal{C}$-free undirected graph $G$ contains a clique or an independent set of size at least 
$|G|^{c_{\mathcal{C}}}$. When $\mathcal{C}$ consists of tournaments we say that $\mathcal{C}$ has the \textit{Erd\H{o}s-Hajnal property} if the following
holds: there exists $c_{\mathcal{C}}>0$ such that every $\mathcal{C}$-free tournament $T$ contains a transitive subtournament of order at least $|T|^{c_{\mathcal{C}}}$.
If $\{H\}$ has the Erd\H{o}s-Hajnal property we slightly violate the notation and simply say that $H$ has the Erd\H{o}s-Hajnal property.

The Erd\H{o}s-Hajnal conjecture states that every undirected graph $H$ or (equivalently) every tournament $H$ has the Erd\H{o}s-Hajnal property. 
One can propose a weaker version of the Erd\H{o}s-Hajnal conjecture, where instead of one graph a class $\mathcal{C}$ of graphs is forbidden.
In particular the following conjecture is open.

\begin{theorem}
\label{erdos-hajnal-conjecture-undirected-relaxed}
For every undirected graph $H$ there exists a constant $c(H)>0$ such that the following holds: every $\{H, H^{c}\}$-free graph $G$ contains a clique or a stable set
of size at least $|G|^{c(H)}$.
\end{theorem}

Several papers mentioned before tackled this problem for specific choices of $H$.
There is a natural corresponding conjecture in the directed setting.

\begin{theorem}
\label{erdos-hajnal-conjecture-directed-relaxed}
For every tournament $H$ there exists a constant $c(H)>0$ such that the following holds: every $\{H, H^{c}\}$-free tournament $T$ contains a transitive subtournament
of size at least $|T|^{c(H)}$.
\end{theorem}

Our paper is the first one that addresses the latter conjecture and proves it for several tournaments $H$ that are not known to have the Erd\H{o}s-Hajnal property.

We need few more definitions.

Let $T$ be a tournament, and let $(v_1, \ldots, v_{|T|})$ be an ordering of its 
vertices; denote this ordering by $\theta$.  We say that an edge $(v_{j},v_{i})$ of $T$ is a \textit{backward edge} under this ordering if $i<j$.  
The \textit{graph of backward edges}  under this ordering,
denoted by $B(T, \theta)$, has vertex set $V(T)$,  and  
$v_i v_j \in E(B(T, \theta))$ if and only if  $(v_i,v_j)$ or 
$(v_j,v_i)$ is a backward edge of $T$ under the ordering $\theta$. 

For an integer $t$, we call the graph $K_{1,t}$ a {\em star}. Let $S$ be a
star with vertex set $\{c, l_1, \ldots, l_t\}$, where $c$ is adjacent to $l_1, \ldots, l_t$. We call $c$ the {\em center of the star}, and
$l_1, \ldots, l_t$ {\em the leaves of the star}. 
Note that in the case $t=1$ we may choose arbitrarily any one of the two vertices to be the center of the star, and the other vertex is then considered to be the leaf. 

A {\em star}
in $B(T, \theta)$ is an induced subgraph with vertex set 
$\{v_{i_0},...,v_{i_{j}},..., v_{i_t}\}$, such that \\
$B(T,\theta)|\{v_{i_0}, \ldots, v_{i_t}\}$ is a star as an undirected graph 
with center $v_{i_j}$, and  $i_0 < \ldots < i_{t}$.  In this case we also  
say that $\{v_{i_0}, \ldots, v_{i_t}\}$ is  a star in $T$.

A {\em right star}
in $B(T, \theta)$ is an induced subgraph with vertex set 
$\{v_{i_0}, \ldots, v_{i_t}\}$, such that \\
$B(T,\theta)|\{v_{i_0}, \ldots, v_{i_t}\}$ is a star with center $v_{i_t}$, 
and  $i_t > i_0, \ldots, i_{t-1}$.  In this case we also  
say that $\{v_{i_0}, \ldots, v_{i_t}\}$ is  a right star in $T$. \\
A {\em left star}
in $B(T, \theta)$ is an induced subgraph with vertex set 
$\{v_{i_0}, \ldots, v_{i_t}\}$, such that 
$B(T,\theta)|\{v_{i_0}, \ldots, v_{i_t}\}$ is a star with center $v_{i_0}$, 
and  $i_0 < i_1, \ldots, i_t$.    In this case we also  
say that $\{v_{i_0}, \ldots, v_{i_t}\}$ is a  left star in $T$.
Finally, a {\em central star}   in $B(T,\theta)$,  is a star that is neither left nor right.

A tournament $T$ is a \textit{galaxy} if there exists an ordering $\theta$ of its vertices such that every connected component of $B(T, \theta)$ is  either a left star, a right stars
or a singleton, and
\begin{itemize}
\item 
no center of a star appears in the ordering between two leaves of another star.
\end{itemize}

We call such an ordering a {\em galaxy ordering} of $T$. 
The following was proved in \cite{choromanski1}.

\begin{theorem}
\label{galaxy-theorem}
Every galaxy satisfies the Erd\H{o}s-Hajnal conjecture.
\end{theorem}

The initial lower bound on the coefficients $\epsilon$ from the conjecture for galaxies  was extremely small (since the initial proof used regularity lemma) but was later significantly improved in a few papers (see \cite{kchoromanski1}, \cite{kchoromanski2}).

The assumption that no center of a star appears in the ordering between two leaves of another star seems artificial
and quite technical yet it is necessary to make proofs work. It is not known whether the conjecture is still true if this
condition is abandoned.
We say that a tournament $T$ is a \textit{nebula} if it has an ordering of vertices $\theta$ such that every connected 
component of $B(T, \theta)$ is a star or a singleton (the star does not have to be necessarily left or right, there is no condition
regarding the location of centers of stars). We call this ordering a {\em nebula ordering} of $T$.   
Notice that every galaxy is obviously a nebula. The following is still a conjecture.

\begin{theorem}
\label{nebula-conjecture}
Every nebula satisfies the Erd\H{o}s-Hajnal conjecture.
\end{theorem}

One may ask whether it helps if it is additionally known that every star of the nebula is very small.
We say that a tournament $T$ is a \textit{left nebula} if it is a nebula and besides under its nebula ordering of vertices
every connected component of the graph of backward edges $B(T,\theta)$ is either a three-vertex left star or a singleton.
We call this nebula ordering a \textit{left nebula ordering}.
Similarly, we say that a tournament $T$ is a \textit{right nebula} if it is a nebula and besides under its nebula ordering of vertices
every connected component of the graph of backward edges $B(T,\theta)$ is either a three-vertex right star  or a singleton.
We call this nebula ordering a \textit{right nebula ordering}.
Finally, we say that a tournament $T$ is a \textit{central nebula} if it is a nebula and besides under its nebula ordering of vertices
every connected component of the graph of backward edges $B(T,\theta)$ is either a three-vertex central star or a singleton.
We call this nebula ordering a \textit{central nebula ordering}.

Unfortunately the following claims are still open.

\begin{theorem}
\label{left-nebula-conjecture}
Every left nebula satisfies the Erd\H{o}s-Hajnal conjecture.
\end{theorem}

\begin{theorem}
\label{right-nebula-conjecture}
Every right nebula satisfies the Erd\H{o}s-Hajnal conjecture.
\end{theorem}

\begin{theorem}
\label{central-nebula-conjecture}
Every central nebula satisfies the Erd\H{o}s-Hajnal conjecture.
\end{theorem}

However if we exclude both: 
\begin{itemize}
\item an arbitrary left nebula and an arbitrary right nebula, or
\item an arbitrary left nebula and an arbitrary central nebula, or
\item an arbitrary right nebula and an arbitrary central nebula,
\end{itemize}

then the conjecture is satisfied.
The main result of this paper states that:

\begin{theorem}
\label{main-theorem}
If $H_{1}$ and $H_{2}$ are:
a left nebula and a right nebula, or:
a left nebula and a central nebula, or:
a right nebula and a central nebula,
then $\{H_{1},H_{2}\}$ has the Erd\H{o}s-Hajnal property.
\end{theorem}

Since one can easily notice that the complement of the left nebula is a right nebula and vice versa, we immediately 
get the following result:

\begin{theorem}
\label{main-theorem-corollary}
If $H$ is a left/right nebula then $\{H,H^{c}\}$ has the Erd\H{o}s-Hajnal property.
\end{theorem}

Tournaments that are not obtained by the mentioned substitution procedure are called $\textit{prime}$.
Tournament that is obtained from two smaller tournaments by the substitution procedure has the
Erd\H{o}s-Hajnal property if these smaller tournaments do have it.
Thus in terms of the conjecture prime tournaments are of main interest.
It can be noticed that the family of left nebulae contains infinitely many tournaments with prime subtournaments $H$ that are neither galaxies nor constellations. For example, take a tournament
$H$ of twelve vertices $\{1,...,12\}$ and with the set of backward edges under ordering $(1,...,12)$
of the form: $\{(5,1),(9,1),(8,6),(11,6),(4,2),(10,3),(12,7)\}$. Using that ordering one can notice that
$H$ is a subtournament of the left nebula. It can be also observed that $H$ is prime and is not a 
constellation nor a galaxy (we leave it to the reader). Similarly, the family of  right nebulae
contains infinitely many tournaments with prime subtournaments $H$ that are neither galaxies
nor constellations. This comes immediately from the previous observation, the fact that the complement
of the left nebula is a right nebula and the fact that the complement of a prime tournament is prime.
Finally, the family of central nebulae contains infinitely many tournaments with prime subtournaments $H$ that are neither galaxies nor constellations. For example, take a tournament
$H$ of twelve vertices $\{1,...,12\}$ and with the set of backward edges under ordering $(1,...,12)$
of the form: $\{(4,1),(8,4),(5,3),(9,5),(6,2),(11,6),(10,7),(12,10)\}$. Using that ordering one can notice that
$H$ is a central nebula. It can be also observed that $H$ is prime and is not a constellation nor a galaxy 
(we leave it to the reader). 
Thus we can conclude that the families of left, right and central nebulae are hard, i.e. are not contained
in the families of tournaments for which the conjecture has been proved so far.
Therefore our results cannot be deduced by the known methods.

This paper is organized as follows:

\begin{itemize}
\item in Section 2 and 3 we present some tools useful in the latter analysis, 
\item in Section 4 we prove Theorem \ref{main-theorem},
\item in the Appendix, for completeness and the convenience of the user, 
        we give the proof of one simple technical result that appeared in a very similar version in another  
        already submitted paper, and that turns out to be useful also in this paper.
\end{itemize}

\section{Product tournaments}

The following notion of a \textit{product tournament} will turn out to be very handy in our further analysis
(we borrow it from \cite{choromanski3}, but for consistency repeat it here).

Let $H_{1},H_{2}$ be two tournaments. Let us consider two injective functions $f_{1}:V(H_{1}) \rightarrow \mathbb{N}, f_{2}:V(H_{2}) \rightarrow \mathbb{N}$.
Assume furthermore that $\forall_{h^{1} \in V(H_{1}),h^{2} \in V(H_{2})}$ we have: $f_{1}(h^{1}) \neq f_{2}(h^{2})$. 
We shortly denote this last condition by: $<f_{1},f_{2}>=0$. Denote by $\theta_{1}$ the ordering of the vertices of $V(H_{1})$ induced by increasing values of $f_{1}$ on  
$V(H_{1})$ and by $\theta_{2}$ the ordering of the vertices of $V(H_{2})$ induced by increasing values of $f_{2}$ on $V(H_{2})$.
Now let us define \textit{the product $H$ of $H_{1}$ and $H_{2}$ under orderings $\theta_{1}$ and $\theta_{2}$} as follows:
\begin{itemize}
\item $V(H)=V(H_{1}) \cup V(H_{2})$,
\item under ordering $\theta$ of $V(H)$ induced by $f_{1},f_{2}$, where: $<f_{1},f_{2}>=0$, the backward edges of $H$ are exactly the backward edges of $H_{1}$ 
      under $\theta_{1}$ and the backward edges of $H_{2}$ under $\theta_{2}$.
\end{itemize}
We denote this product tournament $H$ by $H^{f_{1}}_{1} \oplus H^{f_{2}}_{2}$.
Notice that the $\oplus$ operation is commutative and associative.

Tournament $H$ is called a \textit{small left star} if it consists of three vertices: $c, l_{1}, l_{2}$ such that under
ordering $(c, l_{1}, l_{2})$ the set of backward edges is of the form $\{(l_{1},c),(l_{2},c)\}$. We call this ordering the
\textit{default ordering of a small left star}.
Similarly, tournament $H$ is called a \textit{small right star} if it consists of three vertices: $l_{1}, l_{2}, c$ such that under
ordering $(l_{1}, l_{2}, c)$ the set of backward edges is of the form $\{(c,l_{1}),(c,l_{2})\}$. We call this ordering the
\textit{default ordering of a small right star}. Finally, tournament $H$ is called a \textit{small central star} if it consists
of three vertices $l_{1}, c, l_{2}$ such that under ordering $(l_{1}, c, l_{2})$ the set of backward edges is of the form 
$\{(c,l_{1}),(l_{2},c)\}$. We call this ordering the \textit{default ordering of a small central star}.

Notice that every left nebula is a subtournament of another left nebula $H$ which is of the form $H=H_{1}^{f_{1}} \oplus ... \oplus H_{r}^{f_{r}}$, where
$r>0$, each $H_{i}$ is a small left star, $<f_{i},f_{j}>=0$ for $i \neq j$ and $f_{i}$ is induced by the default ordering of $H_{i}$.
Similarly, every right nebula is a subtournament of another right nebula $H$ which is of the form $H=H_{1}^{f_{1}} \oplus ... \oplus H_{r}^{f_{r}}$, where
$r>0$, each $H_{i}$ is a small right star, $<f_{i},f_{j}>=0$ for $i \neq j$ and $f_{i}$ is induced by the default ordering of $H_{i}$.
Finally, every central nebula is a subtournament of another central nebula $H$ which is of the form $H=H_{1}^{f_{1}} \oplus ... \oplus H_{r}^{f_{r}}$, where
$r>0$, each $H_{i}$ is a small central star, $<f_{i},f_{j}>=0$ for $i \neq j$ and $f_{i}$ is induced by the default ordering of $H_{i}$.
Thus it suffices to prove our main result only for left/right/central nebulae of the form above.
From now one our analysis regards only these types of nebulae.

\section{$(c, \lambda, w)$-{\em structures}}

For a tournament $T$ and two disjoint nonempty sets $A,B \subseteq(V(T))$ let $d(A,B) = \frac{e(A,B)}{|A||B|}$, where $e(A,B)$ is a number of edges from $A$ to $B$.
By $tr(T)$ we denote the size of the largest transitive subtournament of $T$.
Let $c>0$, $0<\lambda<1$ be constants, and let $w$ be a $\{0,1\}$-vector of length $|w|$. Let $T$ be a tournament with $|T|=n$. A sequence of disjoint subsets 
$(S_{1},S_{2},...,S_{|w|})$ of $V(T)$ is a $(c, \lambda, w)$-{\em structure} if
\begin{itemize} 
\item whenever $w_i=0$ we have $|S_{i}| \geq cn$ 
\item whenever $w_i=1$ the set $T|S_{i}$ is transitive and $|S_{i}| \geq c \cdot tr(T)$
\item $d(S_{i},S_{j}) \geq 1 - \lambda$ for all $1 \leq i < j \leq |w|$.
\end{itemize} 

We say that a $(c, \lambda, w)$-{\em structure} $(S_{1},S_{2},...,S_{|w|})$ is \textit{strong} if the following holds for every $v \in S_{i}$, $i=1,...,|w|$ and $j \neq i$:

\begin{itemize}
\item $d(\{v\},S_{j}) \geq 1 - \lambda$ if $i<j$,
\item $d(S_{j},\{v\}) \geq 1 - \lambda$ if $i>j$.
\end{itemize}

In this paper we will only use $(c,\lambda,w)$-structures for vectors $w$ with all entries equal to $0$. 
We shortly denote them as \textit{$(c,\lambda)$-structures}.
However since this construction was first
defined in \cite{choromanski1}, we gave here the most general definition.

\section{Excluding two nebulae}

Our main goal of this section is to prove Theorem \ref{main-theorem}. Before doing it we will introduce
few more definitions and useful technical lemmas.

Let $T$ be a tournament.
Let $\sigma=(S_{1},...,S_{k})$ be a sequence of pairwise disjoint subsets of $V(T)$. 
For $v \in S_{i}$ and $i \neq j$ denote:

\[ N^{\sigma}(v,j) = \left\{
  \begin{array}{l l}
    \{w \in S_{j} : (w,v) \in E(T)\} & \quad \text{if $j > i$},\\
     \{w \in S_{j} : (v,w) \in E(T)\} & \quad \text{if $j < i$}
  \end{array} \right.\]

Let $\sigma = (S_{1},S_{2},S_{3})$ be an ordered triple of pairwise disjoint subsets of $V(T)$.
We say that $\sigma$ is a $(i,j)$-triple, where $i,j \in \{1,2,3\}$, $i \neq j$ if the following holds:

there exists an ordering of the vertices of $S_{i}$: $v_{1}^{i},...,v^{i}_{|S_{i}|}$ such that:
$$\min \{k: (N^{\sigma}(v_{1}^{i},j)+...+N^{\sigma}(v_{k}^{i},j) \geq \frac{|S_{j}|}{2}) \}
\geq \min \{k: (N^{\sigma}(v_{1}^{i},l)+...+N^{\sigma}(v_{k}^{i},l) \geq \frac{|S_{l}|}{2}) \}$$
for $l = 6 - i - j$ and both minima are finite.

Let $\chi=(S_{1},...,S_{k})$ be a strong $(c,\lambda)$-structure.
Assume furthermore that $|S_{1}|=...=|S_{k}|=t$.
Let $H$ be a tournament with $V(H)=\{h_{1},...,h_{|H|}\}$. Assume that there exist indices
$1 \leq i_{1} < ... < i_{|H|} \leq k$ and a bijection $\phi : V(H) \rightarrow \{i_{1},...,i_{|H|}\}$
such that the following holds:

\begin{itemize}
 \item for $h=h_{1},...,h_{|H|}$ there exists an ordering of the vertices $(v_{1}^{h},...,v_{t}^{h})$ of $S_{\phi(h)}$, where $v_{j}^{h} \in S_{\phi(h)}$ for $j=1,...,t$ and
 \item for every $j \in \{1,...,t\}$ a set $\{v_{j}^{h_{1}},...,v_{j}^{h_{|H|}}\}$ induces a copy of $H$ where the isomorphism
       is given by the mapping: $h_{i} \rightarrow v_{j}^{h_{i}}$ for $i=1,...,|H|$.
\end{itemize}

Then we say that $\chi$ is \textit{$(H,\phi)$-normal}.
We start with our first technical lemma.

\begin{theorem}
\label{first-technical-lemma}
Let $\chi$ be a strong $(c,\lambda)$-structure. Assume furthermore that $|S_{1}|=...=|S_{k}|=t$.
Let $H_{1},...,H_{p}$ be tournaments with $V(H_{j})=\{h^{j}_{1},...,h^{j}_{|H_{j}|}\}$ for $j=1,...,p$.
Assume that $\chi$ is $(H_{j},\phi_{j})$-normal for $j=1,...,p$ and $\phi_{j_{1}}(h^{j_{1}}_{i_{1}}) \neq \phi_{j_{2}}(h^{j_{2}}_{i_{2}})$
for $j_{1} \neq j_{2}$, $i_{1} \in \{1,...,|H_{j_{1}}|\}$, $i_{2} \in \{1,...,|H_{j_{2}}|\}$.
Then either $\lambda \geq \frac{1}{(p-1)^{2}(\max_{j=1,...,p} |H_{j}|)^{2}}$ or the tournament induced by 
$S_{1} \cup ... \cup S_{k}$ has a tournament $H = H_{1}^{\phi_{1}} \oplus ... \oplus H_{p}^{\phi_{p}}$ as an induced subtournament.
\end{theorem}

\Proof
Fix some $j$.
We can find $t$ sets: 
$\mathcal{L}^{j}_{1} = \{v_{1}^{j,h_{1}},...,v_{1}^{j, h_{|H_{j}|}}\}$,...,$\mathcal{L}^{j}_{t} = \{v_{t}^{j,h_{1}},...,v_{t}^{j, h_{|H_{j}|}}\}$
such that: tournaments $R_{1}^{j},...,R_{t}^{j}$ induced by $\mathcal{L}^{j}_{1}$,...,$\mathcal{L}^{j}_{t}$ respectively are
isomorphic to $H_{j}$, the isomorphism for any given $\mathcal{L}^{j}_{s}$ is given by the mapping: $h_{i}^{j} \rightarrow v_{s}^{j,h_{i}}$ and 
$v_{s}^{j,h_{i}} \in S_{\phi_{j}(h_{i})}$.
The above holds since $\chi$ is $(H_{j},\phi_{j})$-normal. 
Now consider an $p$-partite undirected graph $B^{p}$ with color classes:

$C_{1}=\{R^{1}_{1},...,R^{1}_{t}\}$,...,$C_{p}=\{R^{p}_{1},...,R^{p}_{t}\}$.
In this graph there exists an edge between vertex $R^{j_{1}}_{s_{1}}$ and $R^{j_{2}}_{s_{2}}$
for $j_{1} < j_{2}$ if the following is true for every
$v_{1} \in V(R^{j_{1}}_{s_{1}})$, $v_{2} \in V(R^{j_{2}}_{s_{2}})$:
\begin{itemize}
\item if $v_{1} \in S_{k_{1}}$, $v_{2} \in S_{k_{2}}$ and $k_{1} < k_{2}$ then $(v_{1},v_{2})$ is an edge,
\item if $v_{1} \in S_{k_{1}}$, $v_{2} \in S_{k_{2}}$ and $k_{1} > k_{2}$ then $(v_{2},v_{1})$ is an edge.
\end{itemize}

Note that since $\chi$ is a strong $(c,\lambda)$-structure, we have for every $j_{1} \in \{1,...,p\}$,
$s \in \{1,...,t\}$ and $j_{2} \neq j_{1}$:
$|N^{j_{2}}(R^{j_{1}}_{s})| \geq |C_{j_{2}}|(1-\lambda |H_{j}|^{2})$, where:
$N^{j_{2}}(R^{j_{1}}_{s})$ is the set of vertices lying in $C_{j_{2}}$ that are neighbors of a vertex $R^{j_{1}}_{s}$ 
in a graph $B^{p}$.
Thus $B^{p}$ has at least $ t^{2}(1-\epsilon){p \choose 2}$ edges, where:
$\epsilon = \lambda (\max_{j}|H_{j}|^{2})$.
Therefore we have: $|E(B^{p})| \geq \frac{p(p-1)}{2} t^{2}(1-\epsilon)=
\frac{(p t)^{2}}{2}(1-\frac{1}{p})(1-\epsilon)=\frac{|V(B^{p})|^{2}}{2}
(1-\frac{1}{p})(1-\epsilon)$.
According to Turan's Theorem (see \cite{diestel}), $B^{p}$ contains a clique of size $p$ if 
$\frac{(1-\frac{1}{p})(1-\epsilon)}{2} > \frac{1}{2}\frac{p-2}{p-1}$, i.e. if $\epsilon < \frac{1}{(p-1)^{2}}$, i.e. if $\lambda < \frac{1}{(p-1)^{2}(\max_{j=1,...,p} |H_{j}|)^{2}}$.
Now note that this clique of size $p$ corresponds to the copy of $H_{1}^{\phi_{1}} \oplus ... \oplus H_{p}^{\phi_{p}}$ in $\chi$.
Indeed, denote this clique as $\mathcal{C}=\{R^{1}_{i_{1}},...,R^{p}_{i_{p}}\}$ for some indices $1 \leq i_{1} < ... < i_{p} \leq p$.
Notice that the set of vertices that is inducing a copy of $H_{1}^{\phi_{1}} \oplus ... \oplus H_{p}^{\phi_{p}}$ is of the form
$\bigcup_{j=1,...,p} V(R^{j}_{i_{j}})$. This observation completes the proof.
\bbox

The following observations will turn out to be crucial for the main proof.

\begin{theorem}
\label{first-simple-technical-lemma}
Let $T$ be a tournament.
Let $\sigma=(S_{1},S_{2},S_{3})$ be a $(2,1)$-triple or a $(3,1)$-triple, where $S_{1},S_{2},S_{3} \subseteq V(T)$.
Assume that $|S_{1}|,|S_{2}|,|S_{3}| \geq cn$ for some $c>0$, where $n=|T|$. Then the following holds:
\begin{itemize}
 \item there exist vertices $v_{1} \in S_{1}, v_{2} \in S_{2}, v_{3} \in S_{3}$ such that
       $(v_{2},v_{1}),(v_{3},v_{1}),(v_{2},v_{3})$ are edges of $T$ or
 \item there exist two disjoint subsets $A,B \subseteq V(T)$ such that $A$ is complete to $B$ and
       $|A|,|B| \geq \frac{c}{2}n$.
\end{itemize}
\end{theorem}

Similarly,

\begin{theorem}
\label{second-simple-technical-lemma}
Let $T$ be a tournament.
Let $\sigma=(S_{1},S_{2},S_{3})$ be a $(2,3)$-triple or a $(1,3)$-triple, where $S_{1},S_{2},S_{3} \subseteq V(T)$.
Assume that $|S_{1}|,|S_{2}|,|S_{3}| \geq cn$ for some $c>0$, where $n=|T|$. Then the following holds:
\begin{itemize}
 \item there exist vertices $v_{1} \in S_{1}, v_{2} \in S_{2}, v_{3} \in S_{3}$ such that
       $(v_{3},v_{1}),(v_{3},v_{2}),(v_{1},v_{2})$ are edges of $T$ or
 \item there exist two disjoint subsets $A,B \subseteq V(T)$ such that $A$ is complete to $B$ and
       $|A|,|B| \geq \frac{c}{2}n$.
\end{itemize}
\end{theorem}

Finally,

\begin{theorem}
\label{third-simple-technical-lemma}
Let $T$ be a tournament.
Let $\sigma=(S_{1},S_{2},S_{3})$ be a $(1,2)$-triple or a $(3,2)$-triple, where $S_{1},S_{2},S_{3} \subseteq V(T)$.
Assume that $|S_{1}|,|S_{2}|,|S_{3}| \geq cn$ for some $c>0$, where $n=|T|$. Then the following holds:
\begin{itemize}
 \item there exist vertices $v_{1} \in S_{1}, v_{2} \in S_{2}, v_{3} \in S_{3}$ such that
       $(v_{2},v_{1}),(v_{3},v_{2}),(v_{1},v_{3})$ are edges of $T$ or
 \item there exist two disjoint subsets $A,B \subseteq V(T)$ such that $A$ is complete to $B$ and
       $|A|,|B| \geq \frac{c}{2}n$.
\end{itemize}
\end{theorem}

Below we prove only \ref{first-simple-technical-lemma}. The proofs of \ref{second-simple-technical-lemma}
and \ref{third-simple-technical-lemma} are completely analogous and we leave it to the reader.

\Proof
Assume first that $\sigma$ is a $(2,1)$-triple. According to the definition of $\sigma$, there exists an ordering of the vertices of
$S_{2}$: $v^{2}_{1},...,v^{2}_{|S_{2}|}$ and some finite $1 \leq k \leq |S_{2}|$ such that:
\begin{equation}
\label{equation_1}
N^{\sigma}(v^{2}_{1},1) \cup ... \cup N^{\sigma}(v^{2}_{k},1) \geq \frac{|S_{1}|}{2}
\end{equation} and
\begin{equation}
\label{equation_2}
N^{\sigma}(v^{2}_{1},3) \cup ... \cup N^{\sigma}(v^{2}_{k},3) \leq \frac{|S_{3}|}{2}.
\end{equation}

Assume that there do not exist vertices $v_{i} \in S_{i}$ for $i=1,2,3$ such that: $(v_{2},v_{1})$,$(v_{3},v_{1})$ and $(v_{2},v_{3})$ are edges.
But then we immediately get that $N^{\sigma}(v^{2}_{r},1)$ is complete to $S_{3} \backslash N^{\sigma}(v^{2}_{r},3)$ for $r=1,...,k$.
Thus we get: $N^{\sigma}(v^{2}_{1},1) \cup ... \cup N^{\sigma}(v^{2}_{k},1)$ is complete to $S_{3} \backslash (N^{\sigma}(v^{2}_{1},3) \cup ... \cup N^{\sigma}(v^{2}_{k},3))$.

Let $A=N^{\sigma}(v^{2}_{1},1) \cup ... \cup N^{\sigma}(v^{2}_{k},1)$, 
$B = S_{3} \backslash  (N^{\sigma}(v^{2}_{1},3) \cup ... \cup N^{\sigma}(v^{2}_{k},3))$.
Then from \ref{equation_1} and \ref{equation_2} we know that $|A| \geq \frac{|S_{1}|}{2}$,
$|B| \geq \frac{|S_{3}|}{2}$. Since each $S_{i}$ is of size at least $cn$ and, as we have noticed before,
$A$ is complete to $B$, we are done for the case when $\sigma$ is a $(2,1)$-triple.
Now assume that $\sigma$ is a $(3,1)$-triple. The proof is similar. According to the definition of $\sigma$, there exists an
ordering of vertices of $S_{3}$: $v^{3}_{1}$,...,$v^{3}_{|S_{3}|}$ and $k$ such that:

\begin{equation}
\label{equation_3}
N^{\sigma}(v^{3}_{1},1) \cup ... \cup N^{\sigma}(v^{3}_{k},1) \geq \frac{|S_{1}|}{2}
\end{equation} and
\begin{equation}
\label{equation_4}
N^{\sigma}(v^{3}_{1},2) \cup ... \cup N^{\sigma}(v^{3}_{k},2) \leq \frac{|S_{2}|}{2}.
\end{equation}

Assume that there do not exist vertices $v_{i} \in S_{i}$ for $i=1,2,3$ such that: $(v_{2},v_{1})$,$(v_{3},v_{1})$ and $(v_{2},v_{3})$ are edges.
But then we immediately get that $N^{\sigma}(v^{3}_{r},1)$ is complete to $S_{2} \backslash N^{\sigma}(v^{3}_{r},2)$ for $r=1,...,k$.
Thus we get: $N^{\sigma}(v^{3}_{1},1) \cup ... \cup N^{\sigma}(v^{3}_{k},1)$ is complete to $S_{2} \backslash (N^{\sigma}(v^{3}_{1},2) \cup ... \cup N^{\sigma}(v^{3}_{k},2))$.

Let $A=N^{\sigma}(v^{3}_{1},1) \cup ... \cup N^{\sigma}(v^{3}_{k},1)$, 
$B = S_{2} \backslash  (N^{\sigma}(v^{3}_{1},2) \cup ... \cup N^{\sigma}(v^{3}_{k},2))$.
Then from \ref{equation_3} and \ref{equation_4} we know that $|A| \geq \frac{|S_{1}|}{2}$,
$|B| \geq \frac{|S_{2}|}{2}$. Since each $S_{i}$ is of size at least $cn$ and, as we have noticed before,
$A$ is complete to $B$, we are done also for the case when $\sigma$ is a $(3,1)$-triple.
That completes the entire proof.

\bbox

The next lemma gives us another scenario where we obtain two disjoint linear sets such that one is complete to the other one.

\begin{theorem}
\label{two-triples-lemma}
Let $T$ be a tournament.
Let $\sigma=(S_{1},S_{2},S_{3})$ be a triple of pairwise disjoint sets of vertices $S_{1},S_{2},S_{3} \subseteq V(T)$. Assume that $|S_{1}|,|S_{2}|,|S_{3}| \geq cn$ for some $c>0$, where $n=|T|$.
Then $\sigma$ is a $(i,j)$-triple or a $(i,6-i-j)$-triple or  there exist two disjoint subsets $A,B \subseteq V(T)$
such that $A$ is complete to $B$ and $|A|,|B| \geq \frac{c}{2}n$.
\end{theorem}

\Proof
Assume that $\sigma$ is not a $(i,j)$-triple and not a $(i,6-i-j)$-triple. Denote $l=6-i-j$.
Then, from the definition of a $(i,k)$-triple we have:
\begin{equation}
\label{equation_5}
\bigcup_{v \in S_{i}} N^{\sigma}(v,j) < \frac{|S_{j}|}{2}
\end{equation}
and
\begin{equation}
\label{equation_6}
\bigcup_{v \in S_{i}} N^{\sigma}(v,l) < \frac{|S_{l}|}{2}.
\end{equation}
But then note that $(S_{j} \backslash \bigcup_{v \in S_{i}} N^{\sigma}(v,j))$ is complete to
$S_{i}$ or $S_{i}$ is complete to $(S_{j} \backslash \bigcup_{v \in S_{i}} N^{\sigma}(v,j))$.
Note also that $|(S_{j} \backslash \bigcup_{v \in S_{i}} N^{\sigma}(v,j))| \geq \frac{|S_{j}|}{2}$
from \ref{equation_5}. Since each $S_{i}$ is of size at least $cn$, we have detected two disjoint sets, each of size at least $\frac{cn}{2}$, such that one is complete to the other one. That completes the proof.
\bbox

Let us define $R^{H}(3,k)$ to be the smallest number such that every $3$-regular hypegraph with 
$n \geq R^{H}(3,k)$ vertices and ${n \choose 3}$ edges, each colored black or white, contains a monochromatic
clique of size at least $k$. By standard Ramsey Theorem, $R^{H}(3,k)$ is finite.

Let $L=L_{1}^{f_{1}} \oplus ... \oplus L_{l}^{f_{l}}$ be a left nebula, where each $L_{i}$ is a small
left star. Let $R=R_{1}^{g_{1}} \oplus ... \oplus R_{r}^{g_{r}}$ be a right nebula, where each $R_{i}$ is a small right star.  Let $C=C_{1}^{h_{1}} \oplus ... \oplus C_{c}^{h_{c}}$ be a central nebula, where each $C_{i}$ is a central star. We have already noticed that it suffices to prove \ref{main-theorem} for nebulae of this form. Without loss of generality we can assume that
$\bigcup_{i=1,...l} f_{i}(V(L_{i})) \subseteq \{1,...,k\}$,
$\bigcup_{i=1,...r} g_{i}(V(R_{i})) \subseteq \{1,...,k\}$
and $\bigcup_{i=1,...c} h_{i}(V(C_{i})) \subseteq \{1,...,k\}$ for some $k>0$.

Below we state the result that will directly lead to \ref{main-theorem}.

\begin{theorem}
\label{main-technical-theorem}
Let $\chi_{0}=(S^{0}_{1},...,S^{0}_{t})$ be a strong $(c,\lambda)$-structure with vertices taken from the $\{L,R\}$-free tournament $T$ and let $t=R^{H}(3,k)$. Assume that $|S^{0}_{1}|=|S^{0}_{2}|=...=|S^{0}_{t}|$. Assume that the following holds:
$$\lambda < \frac{1}{3k{R^{H}(3,k) \choose k}(l-1)^{2}(r-1)^{2}(\max_{j=1,...,l}|L_{j}|^{2})(\max_{j=1,...,r}|R_{j}|)^{2}}.$$
Then there exist two disjoint subsets $A,B \subseteq V(T)$ such that $|A|,|B| \geq \frac{c}{6}n-3k{t \choose k}$ and
$A$ is complete to $B$.
Similar result holds if we replace $\{L,R\}$-free tournament $T$ by $\{L,C\}$-free tournament $T$
(then instead of $R_{j}$ and $r$ we have $C_{j}$ and $c$ in the upper bound on $\lambda$) or by $\{R,C\}$-free tournament $T$ (then instead of $L_{j}$ and $l$ we have $C_{j}$ and $c$ in the upper bound on $\lambda$).

\end{theorem}

Before proving this technical theorem we will show why it leads to the proof of \ref{main-theorem}.
Let us assume that \ref{main-technical-theorem} holds and conduct the proof of \ref{main-theorem}.

\Proof
Note first that for any $p \in \mathbb{N}$, any $\lambda > 0$ and any tournament $H$ there exist:
$c>0$ and $n_{0} > 0$ such that every $H$-free tournament $T$ with $|T| \geq n_{0}$ vertices contains
a strong $(c,\lambda)$-structure $\chi=(S_{1},...,S_{p})$, where $|S_{1}|=...=|S_{p}|$. Very similar 
observation appeared in \cite{choromanski1}, but for completeness we will give the entire proof of that
fact in the Appendix. We will take this fact for granted now.
Without loss of generality we can assume that the following holds:
$n_{0} \geq \frac{36k{R^{H}(3,k) \choose k}}{c}$, where $k$ was defined just before the statement of \ref{main-technical-theorem} was given.
Let us denote $\epsilon = \min(\log_{n_{0}}(2),\log_{\frac{c}{12}}(\frac{1}{2}))$. Notice first that trivially every tournament $T$ with less than
$n_{0}$ vertices contains a transitive subtournament of size at least $|T|^{\epsilon}$.
Now let $T$ be a $\{H_{1},H_{2}\}$-free tournament with  $n \geq n_{0}$ vertices, where either:
$H_{1}$ is a left nebula and $H_{2}$ is a right nebula or: $H_{1}$ is a left nebula and $H_{2}$ is a central
nebula or: $H_{1}$ is a right nebula and $H_{2}$ is a central nebula.
We want to prove that $T$ contains a transitive subtournament of order at least $|T|^{\epsilon}$.
We will proceed by induction on $|T|$. The base case for $n \leq n_{0}$ was already proven. Now assume that for every $\{H_{1},H_{2}\}$-free tournament on less than $n$ vertices the statement is true.
Take $n$-vertex tournament $T$. By our earlier observation, $T$ contains a strong $(c,\lambda)$-structure
$\chi=(S_{1},...,S_{p})$ for $p=R^{H}(3,k)$ and $\lambda = \frac{1}{6k{R^{H}(3,k) \choose k}(l-1)^{2}(r-1)^{2}(\max_{j=1,...,l}|L_{j}|^{2})(\max_{j=1,...,r}|R_{j}|)^{2}}$. But then, according to \ref{main-technical-theorem},
there exist two disjoint sets $A,B \subseteq V(T)$ such that $A$ is complete to $B$ and 
$|A|,|B| \geq \frac{c}{6}n-3k{R^{H}(3,k) \choose k} \geq \frac{cn}{12}$. By induction: $tr(T|A) \geq |A|^{\epsilon} \geq (\frac{c}{12}n)^{\epsilon}$ and $tr(T|B) \geq |B|^{\epsilon} \geq (\frac{c}{12}n)^{\epsilon}$ (notice that from the choice of $n_{0}$ we have: $ 3k{R^{H}(3,k) \choose k} \leq \frac{cn}{12}$).
But then combining the largest transitive subtournament of $T|A$ with the largest transitive subtournament of $T|B$ we obtain transitive subtournament $W$ of order at least $2(\frac{c}{12}n)^{\epsilon}$. 
Thus, from the choice of $\epsilon$, we get: $|W| \geq |T|^{\epsilon}$. That completes the proof.
\bbox

Thus it remains to prove \ref{main-technical-theorem}.

\Proof
Let us denote $W=|S^{0}_{1}|=...=|S^{0}_{t}|$.
We will show the proof for a $\{L,R\}$-free tournament $T$. For the two remaining cases the analysis
is completely analogous.
We need to introduce few more useful objects. Let $\mathcal{C}$ be the set of all ${t \choose k}$
subsets of the set $\{1,...,t\}$ of size $k$.
Denote $\mathcal{C}=\{C^{1},...,C^{{t \choose k}}\}$ and $C^{i}=\{c^{i}_{1},...,c^{i}_{k}\}$,
where $c^{i}_{1} < ... < c^{i}_{k}$.
We will construct subsets $A,B$ from the statement of the theorem algorithmically.
Our algorithm consists of several phases. In each phase for every $k$-element subset $C^{i}$
we keep two vectors: $v_{R}^{i}$ and $v_{L}^{i}$. These vectors change over time during the execution
of the algorithm. We will explain later how. In every phase of the algorithm $v_{R}^{i}$ is of length $r$ and 
$v_{L}^{i}$ is of length $l$. Throughout the execution of the algorithm each entry of $v_{R}^{i}$
and each entry of $v_{L}^{i}$ is a set of pairwise disjoint $3$-element ordered triples of vertices of $\chi_{0}$ (disjointness is in the set sense, i.e. $3$-elements sets corresponding to $3$-element triples
after forgetting the order are pairwise disjoint).
Let $v_{R}^{i}(j)$ be the $j^{th}$ element of $v_{R}^{i}$ and let $v_{L}^{i}(j)$ be the $j^{th}$ element of $v_{L}^{i}$. Let $V(L_{j})=\{l^{j}_{1},l^{j}_{2},l^{j}_{3}\}$, where $(l^{j}_{1},l^{j}_{2},l^{j}_{3})$
is a default ordering of a small left star. Each element of $v_{L}^{i}(j)$ (notice that $v_{L}^{i}(j)$ is itself a set) is of the form
$(w_{1},w_{2},w_{3})$, where: $\{w_{1},w_{2},w_{3}\}$ induces a tournament isomorphic to $L_{j}$, $(w_{1},w_{2},w_{3})$ is a default ordering of a small left star and and besides: $w_{i} \in S^{0}_{q}$ for $q=c^{i}_{f_{j}(l_{i}^{j})}$.
Let $V(R_{j})=\{r^{j}_{1},r^{j}_{2},r^{j}_{3}\}$, where $(r^{j}_{1},r^{j}_{2},r^{j}_{3})$
is a default ordering of a small right star. Each element of $v_{R}^{i}(j)$ is of the form
$(z_{1},z_{2},z_{3})$, where: $\{z_{1},z_{2},z_{3}\}$ induces a tournament isomorphic to $R_{j}$, $(z_{1},z_{2},z_{3})$ is a default ordering of a small right star and and besides: $z_{i} \in S^{0}_{u}$ for $u=c^{i}_{g_{j}(r_{i}^{j})}$.
The above properties of vectors $v_{R}^{i}$,$v_{L}^{i}$ will be valid during the entire execution of the algorithm. 
Furthermore, during the entire execution of the algorithm each entry of $v_{R}^{i}$
and each entry of $v_{L}^{i}$ will be of size at most $\lceil \frac{W}{9k{t \choose k}} \rceil$.
When we describe our algorithm in details it will be easy to check that all the properties above are satisfied. 
Before going into details of the algorithm let us notice one fundamental fact - under the assumptions above it is not possible that at same phase of the algorithm there exists some $i$ such that each entry of $v_{R}^{i}$ is of size 
$\lceil \frac{W}{9k{t \choose k}} \rceil$ or each entry of $v_{L}^{i}$ is of size $\lceil \frac{W}{9k {t \choose k}} \rceil$.
We call this property the \textit{nonsaturation property}. Let us understand why the nonsaturation property holds. Assume by contradiction that it does not hold, i.e. for some $i$ all entries of $v_{R}^{i}$
are of size $\lceil \frac{W}{9k {t \choose k}} \rceil $ or all entries of $v_{L}^{i}$ are of size at $\lceil \frac{W}{9k {t \choose k}} \rceil$.
Without loss of generality assume that the latter is true. 
Let $\mathcal{L}=\bigcup_{j} \{v : \exists_{x,y} (v,x,y) \in v_{L}^{i}(j)\} \cup \bigcup_{j} \{v : \exists_{x,y} (x,v,y) \in v_{L}^{i}(j)\} \cup \bigcup_{j} \{v : \exists_{x,y} (x,y,v) \in v_{L}^{i}(j)\}$. Denote $I_{j} = \mathcal{L} \bigcap S_{c^{i}_{j}}$ for $j=1,...k$.
Let $\{I_{u_{1}},...,I_{u_{m}}\}$ be those $I_{j}$ that are nonempty ($u_{1} < ... < u_{m}$).
Since the nonsaturation property is not satisfied, we get: $|I_{u_{1}}|,...,|I_{u_{m}}| = \lceil \frac{W}{9k{t \choose k}} \rceil$. 
Note that $(I_{u_{1}},...,I_{u_{m}})$ is a strong $(c \theta,\frac{\lambda}{\theta})$-structure, where: $\theta=\frac{1}{3k{t \choose k}}$. This comes directly from derived lower bound
on the sizes of $I_{u_{j}}$ and the fact that $\chi_{0}$ is a strong $(c,\lambda)$-structure.
Denote this $(c \theta,\frac{\lambda}{\theta})$-structure by $\Omega=(S^{\Omega}_{1},...,
S^{\Omega}_{m})$.
Notice, that from the definition of  $(S^{\Omega}_{1},...,
S^{\Omega}_{m})$ and $v_{L}^{i}$
we know that $\Omega$ is $(L_{j}, \xi_{j})$-normal for $j=1,...,l$ for some $\xi_{j}$. 
Indeed, $\xi_{j}$ is defined as follows: $\xi_{j}(l^{j}_{i})=y$ such that $u_{y} = c^{i}_{f_{j}(l^{j}_{i})}$. But then, since $\frac{\lambda}{\theta} < \frac{1}{(l-1)^{2}}(\max_{j=1,...,l} |L_{j}|)^{2}$, by \ref{first-technical-lemma} we get: $\Omega$ contains a copy
of $L^{f_{1}}_{1} \oplus ... \oplus L^{f_{l}}_{l}$ which is a contradiction.
We obtain similar contradiction if we assume that every entry of $v_{R}^{i}$ is of size
$\lceil \frac{W}{9k {t \choose k}} \rceil$.

Thus we can conclude that the nonsaturation property is satisfied.

Now we will finally describe the algorithm in details.
In our algorithm we will use the $3$-regular hypergraph $\mathcal{H}$ with the set
of vertices $V(\mathcal{H})=\{1,...,t\}$ and all possible ${t \choose k}$ egdes.  
The algorithm maintains the sequence $(S_{1},...,S_{t})$ and modifies it throughout its execution. The initial $(S_{1},...,S_{t})$ is of the form: $(S^{0}_{1},...,S^{0}_{t})$,
the sequence before the $i^{th}$ phase of the algorithm will be denoted as
$(S^{i}_{1},...,S^{i}_{t})$ for $i=0,1,...$.
The $i^{th}$ phase of the algorithm is conducted as follows. 
For every edge $e=\{j_{1},j_{2},j_{3}\}$ of $\mathcal{H}$, where $j_{i} < j_{2} < j_{3}$ we color
it as white if $(S^{i}_{j_{1}},S^{i}_{j_{2}},S^{i}_{j_{3}})$ is a $(2,1)$-triple
and as black if it is a $(2,3)$-triple. Note that from \ref{two-triples-lemma} we know that
if $e$ is not colored then two disjoint sets $A$, $B$, each of size at least $\frac{\min_{k=1,2,3}|S^{i}_{j_{k}}|}{2}$ such that $A$ is complete to $B$, were detected.
If this is the case then we terminate the algorithm and say that \textit{state $0$ was reached}.
Assume therefore that this is not the case. In such a scenario every edge of $\mathcal{H}$ is colored either black or white. Since $\mathcal{H}$ has $R^{H}(3,k)$ vertices, from the definition
of $R^{H}$ we conclude that it contains a monochromatic clique $C_{j}=\{c^{j}_{1},...,c^{j}_{k}\}$
for some $1 \leq c^{j}_{1} < ... < c^{j}_{k} \leq t$. Without loss of generality we can assume that all vertices of $C_{j}$ are white. Now let us take vector $v_{L}^{j}$.
If each entry of $v_{L}^{j}$ is of size at least $\lceil \frac{W}{9k{t \choose k}} \rceil$ then we
terminate our algorithm and say that the algorithm \textit{reached state 1}. 
Otherwise we find an arbitrary entry $z$ of $v_{L}^{j}$ of size smaller than $\lceil \frac{W}{9k{t \choose k}} \rceil$. Assume that it corresponds to the left star $L_{k}$ with 
$V(L_{k}) = \{l^{k}_{1},l^{k}_{2},l^{k}_{3}\}$, where $(l^{k}_{1},l^{k}_{2},l^{k}_{3})$
is a default ordering of a small left star.
Let us now take triple $(S^{i}_{x_{1}},S^{i}_{x_{2}},S^{i}_{x_{3}})$, where:
$x_{1}=c^{j}_{f_{k}(l^{k}_{1})}$, $x_{2}=c^{j}_{f_{k}(l^{k}_{2})}$,
$x_{3}=c^{j}_{f_{k}(l^{k}_{3})}$. This is a $(2,1)$-triple since an edge $\{x_{1},x_{2},x_{3}\}$ is white. Then, using \ref{first-simple-technical-lemma}, we can conclude that we either find
two disjoint subsets $A,B$ such that $A$ is complete to $B$ and 
$|A|,|B| \geq \frac{\min_{k=1,2,3} |S_{x_{k}}|}{2}$ or we find three vertices $v_{1},v_{2},v_{3}$, as in the statement of \ref{first-simple-technical-lemma}.
In the former case the algorithm is terminated and we say that it \textit{reached state 2}.
In the latter case triple $(v_{1},v_{2},v_{3})$ is added to $v^{j}_{L}(z)$ and the sequence
$(S^{i+1}_{1},...,S^{i+1}_{t})$ is obtained from $(S^{i}_{1},...,S^{i}_{t})$ by deleting
vertices $v_{1},v_{2},v_{3}$ from those sets of $\{S^{i}_{1},...,S^{i}_{t}\}$ that contain them.
If all edges of $C^{j}$ are black then the analysis is completely analogous but instead of a small left star and $(2,1)$-triple definition we use small right star and $(2,3)$-triple definition. We again reach state 1 or 2
or delete from $\Omega$ vertices $v_{1},v_{2},v_{3}$ inducing small right star. In the latter scenario, similarly as previously, $(v_{1},v_{2},v_{3})$ is added to $v^{j}_{R}(z)$.

Let us analyze the algorithm. First notice that all the properties of the algorithm mentioned earlier (in particular those regarding vectors $v_{L}^{i}$ and $v_{R}^{i}$) are trivially satisfied. Thus, as we have observed earlier, the nonsaturation property holds. But that implies in particular that the algorithm cannot reach state $1$. Note also that the algorithm has to terminate.
Indeed, this follows directly from: 
\begin{itemize}
 \item the nonsaturation property, 
 \item the observation that at each phase of the algorithm as long as it does not terminate
         at least one entry of some $v_{R}^{i}$ or some $v_{L}^{i}$ increases in size and
 \item the fact that throughout the execution of the algorithm each entry of each $v_{R}^{i}$
         and $v_{L}^{i}$ is of size at most $\lceil \frac{W}{9k {t \choose k}} \rceil$.
\end{itemize}
Thus the algorithm terminates at state 0 or 2. In both cases two disjoint subsets $A,B$ such that
$A$ is complete to $B$ are detected. Let us find upper bound on their sizes. 
To do that let us fix a set $S^{0}_{i}$ for some $1 \leq i \leq t$.
Throughout the execution of the algorithm some vertices may be deleted from $S^{0}_{i}$
and be moved to entries of the vectors $v_{L}^{i}$ and $v_{R}^{i}$. Notice however that the total number of these entries is at most $2k{t \choose k}$ (since the length of each $v_{L}^{i}$ and each $v_{R}^{i}$ is at most $k$ and $i \in \{1,...,{t \choose k}\}$).
Besides, as noticed before, each entry is of size at most $\lceil \frac{W}{9k {t \choose k}} \rceil$.
Thus the total number of elements deleted from $S^{0}_{i}$ throughout the execution of the
algorithm is at most: $3 \cdot 2k{t \choose k} \lceil \frac{W}{9k {t \choose k}} \rceil \leq 
\frac{2}{3}W + 6k{t \choose k}$. Thus throughout the execution of the algorithm
each $S^{i}_{j}$ is of size at least $\frac{W}{3} - 6k{t \choose k}$.
Thus in particular at state 0 or 2 the subsets $A$ and $B$ that are found are of size at least
$\frac{W}{6} - 3k{t \choose k} \geq \frac{cn}{6}-3k{t \choose k}$ each. That completes the proof.

In the scenario where we exclude the left nebula $L$ and a central nebula $C$ the proof is almost exactly the same. Instead of \ref{second-simple-technical-lemma} we are using \ref{third-simple-technical-lemma}. Instead of $(2,1)$-triples
and $(2,3)$-triples we are using $(3,1)$-triples and $(3,2)$-triples.
In the scenario where we exclude the right nebula $R$ and the central nebula $C$ the proof is also almost exactly the same. The difference now is that we are using \ref{third-simple-technical-lemma} and besides $(1,2)$- and $(1,3)$- triples.
\bbox

\section{Appendix}

The goal of this section is to prove the following result used in the main body of the paper:

\begin{theorem}
\label{simplelemma}
Let $H$ be a tournament, $P>0$ be an integer and $0 < \lambda < 1$. Then there is an integer $N$ such that for every tournament $T$ not containing $H$  
and with $|T| \geq N$ there exists a constant $c>0$ and $P$ pairwise disjoint subsets $A_{1},A_{2},...,A_{P}$ 
of the vertices of $T$ satisfying:
\begin{itemize}
\item $d(\{v\},A_{j}) \geq 1-\lambda$ for $i,j \in \{1,2,...,P\}$, $i<j$, $v \in A_{i}$,
\item $d(A_{i},\{v\}) \geq 1-\lambda$ for $i,j \in \{1,2,...,P\}$, $i<j$, $v \in A_{j}$,
\item $|A_{i}| \geq c|T|$ for $i \in \{1,2,...,P\}$,
\item $|A_{1}|=...=|A_{P}|$.
\end{itemize}
\end{theorem}

Before accomplishing it we need to introduce few more definitions.

Let $T$ be a tournament.
Given $\epsilon>0$ we call a pair ${A,B}$ of disjoint subsets of $V(T)$ 
$\epsilon$-{\em regular} if all $X \subseteq A$ and $Y \subseteq B$ with $|X| \geq \epsilon |A|$ and $|Y| \geq \epsilon |B|$ satisfy:
$|d(X,Y)-d(A,B)| \leq \epsilon$.

Consider a partition $\{V_{0},V_{1},...,V_{k}\}$ of $V(T)$ in which one set $V_{0}$ has been singled out as an \textit{exceptional} set. (This exceptional set $V_{0}$ may be empty). We call such a partition an $\epsilon$-regular partition of $T$ if it satisfies the following three conditions: 
\begin{itemize}
\item $|V_{0}| \leq \epsilon |V|$ 
\item $|V_{1}|=...=|V_{k}|$
\item all but at most $\epsilon k^{2}$ of the pairs $(V_{i},V_{j})$ with $1 \leq i < j \leq k$ are $\epsilon$-regular.
\end{itemize}

The following was proved in \cite{alon}:

\begin{theorem}
\label{regularitylemma}
For every $\epsilon>0$ and every $m \geq 1$ there exists an integer $DM=DM(m,\epsilon)$ such that every tournament of order at least $m$ admits an $\epsilon$-regular partition $\{V_{0},V_{1},...,V_{k}\}$ with $m \leq k \leq DM$.
\end{theorem}

The above lemma is a ``tournament"-version of the celebrated \textit{Regularity Lemma} proved by Endre Szem\'{e}redi and originally stated for undirected graphs
(\cite{diestel2}). In the undirected setting we only need to change the definition of $e_{X,Y}$ which is now the number of edges between sets $X$ and $Y$.
The original version of the lemma is as follows:

\begin{theorem}
\label{regularitylemma_original}
For every $\epsilon>0$ and every $m \geq 1$ there exists an integer $DM=DM(m,\epsilon)$ such that 
every undirected graph of order at least $m$ admits an $\epsilon$-regular partition $\{V_{0},V_{1},...,V_{k}\}$ with $m \leq k \leq DM$.
\end{theorem}

We also need the following lemma:

\begin{theorem}
\label{universalgraphlemma}
For every natural number $k$ and real number $0< \lambda < 1$ there exists $0 < \eta = \eta(k, \lambda) < 1$ such that for every tournament $H$ 
 with vertex set $\{x_{1},...,x_{k}\}$ and tournament 
$T$ with vertex set $V(T)=\bigcup_{i=1}^{k} V_{i}$, if the $V_{i}$'s are disjoint sets, each of order at least one, and each pair $(V_{i},V_{j})$, $1 \leq i < j \leq k$ is $\eta$-regular, with $d(V_{i},V_{j}) \geq \lambda$ and $d(V_{j},V_{i}) \geq \lambda$, then there exist vertices $v_{i} \in V_{i}$ for $i \in \{1, \ldots, k\}$,  
such that the map $x_{i} \rightarrow v_{i}$ gives an isomorphism between $H$ and the subtournament of $T$ induced by $\{v_{1},...,v_{k}\}$.
\end{theorem}

The undirected version of the lemma above is another celebrated result, the so-called \textit{Embedding Lemma}. 

\begin{theorem}
\label{universalgraphlemma_original}
For every natural number $k$ and real number $0< \lambda < 1$ there exists $0 < \eta = \eta(k, \lambda) < 1$ such that for every undirected graph $H$ 
 with vertex set $\{x_{1},...,x_{k}\}$ and undirected graph 
$T$ with vertex set $V(T)=\bigcup_{i=1}^{k} V_{i}$, if the $V_{i}$'s are disjoint sets, each of order at least one, and each pair $(V_{i},V_{j})$, 
$1 \leq i < j \leq k$ is $\eta$-regular, with $d(V_{i},V_{j}) \geq \lambda$ and $d(V_{j},V_{i}) \geq \lambda$, then there exist vertices $v_{i} \in V_{i}$ 
for $i \in \{1, \ldots, k\}$,  
such that the map $x_{i} \rightarrow v_{i}$ gives an isomorphism between $H$ and the subgraph of $T$ induced by $\{v_{1},...,v_{k}\}$.
\end{theorem}

Its proof can be found in \cite{bollobas}.
We will omit the proof of \ref{universalgraphlemma} since it is completely analogous to the proof of the Embedding Lemma.

We are ready to prove \ref{simplelemma}.

\Proof
Write $|T|=n$, $|H|=h$.
Let $R(t_{1},t_{2})$ denote the smallest integer such that every graph of order at least $R(t_{1},t_{2})$ contains either a stable set of size $t_{1}$ or a clique of size $t_{2}$ (by Ramsey theory, $R(t_{1},t_{2})$ is finite). Take $k=R(2^{P-1},h)$.
Take $\eta=\min(\frac{1}{2(k-1)},\eta_{0}(h,\Lambda))$ for $\Lambda = \frac{\lambda}{4P}$ (where $\eta_{0}$ is as in the statement of \ref{universalgraphlemma}).
Let $u>0$ be the smallest integer such that: ${\hat{u} \choose 2} - \eta \hat{u}^{2} > \frac{1}{2}\frac{k-2}{k-1}\hat{u}^{2}$ holds for all $\hat{u} \geq u$. 
By \ref{regularitylemma} there exists an integer $N>0$ such that  every tournament $T$ with $|T| \geq N$ admits an $\eta$-regular partition with at least $u$ parts. Denote by $DM$ the upper bound (from \ref{regularitylemma}) on the number of parts of this partition. Denote the parts of the partition by: $W_{0}, W_{1},...,W_{r}$, where $u \leq r \leq DM$ and $W_{0}$ is the exceptional set.
We have: $|W_{i}| \geq \frac{(1-\eta)n}{DM}$ and besides: $|W_{1}|=...=|W_{r}|$.
Now consider the graph $G$ with $V(G)=\{W_1, \ldots, W_r\}$ where
there is an edge between two vertices if the pair  $(W_{i},W_{j})$ is $\eta$-regular. Then, from the definition of $u$, we have: $|E(G)| \geq \frac{k-2}{2(k-1)}|V(G)|^{2}$. So by Turan's theorem (see \cite{diestel}) it follows that $G$ has a clique of size at least $k$. That means that there exist $k$ parts of the partition, without loss of generality $W_{1},...W_{k}$, such that for all $i,j \in \{1,2,...,k\}, i \neq j$ the pair $(W_{i},W_{j})$  is $\eta$-regular. 
We say that a pair $(W_{i},W_{j})$ for $i,j \in \{1,2,...,k\}, i \neq j$ is \textit{good} if $\Lambda \leq d(W_{i},W_{j}) \leq 1-\Lambda$. Otherwise we say this pair is \textit{bad}.
Now consider the graph $\hat{G}$ with $V(\hat{G})=\{W_{1},...W_{k}\}$, where there is an edge between $W_{i}$ and $W_{j}$ for $i,j \in \{1,...,k\}$, $i \neq j$ if $(W_{i},W_{j})$ is a good pair.
From the definition of $k$ we know that $\hat{G}$ contains a clique of size $h$ or a stable set of size $2^{P-1}$. In other words, either
\begin{itemize}
\item there exist $h$ parts of the partition, without loss of generality denote them $W_{1},...W_{h}$ such that every pair $(W_{i},W_{j})$ is $\eta$-regular and $\Lambda \leq d(W_{i},W_{j}) \leq 1 - \Lambda$ for $i,j \in \{1,2,...,k\}, i \neq j$, or
\item there exist $2^{P-1}$ parts of the partition, without loss of generality denote them $W_{1},...W_{2^{P-1}}$ such that every pair $(W_{i},W_{j})$ is $\eta$-regular and $d(W_{i},W_{j}) > 1 - \Lambda$ or $d(W_{j},W_{i}) > 1  -\Lambda$  for $i,j \in \{1,2,...,2^{P-1}\}, i \neq j$.
\end{itemize}

Since $T$ does not contain $H$ and $\eta \leq \eta_{0}$, \ref{universalgraphlemma} implies that the former is impossible.

Now define $\hat{T}$ to be the tournament with  $V(\hat{T})=\{W_{1},...,W_{2^{P-1}}\}$, where an edge is directed from $W_{i}$ to $W_{j}$ if $d(W_{i},W_{j}) > 1 - \Lambda$ and from $W_{j}$ to $W_{i}$ otherwise.
Using the fact that every tournament of order at least $2^{P-1}$ contains a transitive subtournament of order at least $P$ (see \cite{stearns}), we conclude that $\hat{T}$ contains 
a 
transitive subtournament of order $P$.
That means that there exist 
$P$ 
parts of the partition, without loss of generality $W_{1},...,W_{P}$, such that $d(W_{i},W_{j}) \geq 1 - \Lambda$ for $i,j \in \{1,2,...,P\}, i < j$.
Note that each $W_{i}$ is of order at least $\frac{(1-\eta)n}{DM}$
and $|W_{1}|=...=|W_{P}|$.
Fix $i \in \{1,...,P\}$. For a given $j \in \{1,...,P\}$, $j \neq i$ define: 
$Q^{i}_{j}$ = $\{v \in W_{i}: d(\{v\},W_{j}) \geq 1 - 2P\Lambda\}$ if $i<j$
and $Q^{i}_{j}$ = $\{v \in W_{i}: d(W_{j},\{v\}) \geq 1 - 2P\Lambda\}$ if $j<i$.
Immediately from the fact that $d(W_{\min(i,j)},W_{\max(i,j)}) > 1 - \Lambda$
and from the definition of $Q^{i}_{j}$, we get: $|Q^{i}_{j}| \geq |W_{i}|(1-\frac{1}{2P})$.
Now for every $i \in \{1,...,P\}$ define $F_{i} = \bigcap_{j \in \{1,...,P\}, j \neq i} Q^{i}_{j}$.
From the lower bound on the size of $Q^{i}_{j}$ we get: $|F_{i}| \geq |W_{i}| -  P \cdot \frac{1}{2P} |W_{i}| \geq \frac{1}{2} |W_{i}|$.
For $i=1,...,P$ let $A_{i}$ be the subset of $\lceil \frac{1}{2} |W_{i}|\rceil$ arbitrarily chosen
elements of $F_{i}$.
Notice that by the definition of $Q^{i}_{j}$ and $F_{i}$ we get the following.
For every $i,j \in \{1,...,P\}$, $i < j$, $v \in A_{i}$ we have $d(\{v\},W_{j}) \geq 1 - 2P\Lambda$.
Similarly, for every $i,j \in \{1,...,P\}$, $i > j$, $v \in A_{i}$ we have: $d(W_{j},\{v\}) \geq 1 - 2P\Lambda$.
Thus, since $|A_{i}| \geq \frac{1}{2} |W_{i}|$, we also immediately get the following.
For every $i,j \in \{1,...,P\}$, $i < j$, $v \in A_{i}$ we have $d(\{v\},A_{j}) \geq 1 - 2 \cdot 2P\Lambda = 1-\lambda$.
Similarly, for every $i,j \in \{1,...,P\}$, $i < j$, $v \in A_{j}$ we have: $d(A_{i},\{v\}) \geq 1 - 2 \cdot 2P\Lambda = 1-\lambda$. Notice also that $|A_{1}|=...=|A_{P}| \geq \frac{(1-\eta)n}{2DM}$.
Thus taking $A_{1},...,A_{P}$ and $c=\frac{(1-\eta)n}{2DM}$ we complete the proof
of \ref{simplelemma}.
\bbox

\end{document}